\documentstyle[12pt]{article}    
\newcounter{props}[section]

\newcommand{\wrt}{with respect to\ }

\newcommand{\tx}{there exist\ }
\newcommand{\txs}{there exists\ }
\newcommand{\st}{such that\ }

\newcommand{\fe}{ for every\  }
\newcommand{\fa}{for each\  }

\newcommand{\pair}{pairwise disjoint }
\newcommand{\sq}{sequence }

\newcommand{\dsq}{decreasing sequence\ }
\newcommand{\lra}{\longrightarrow }

\newcommand{\cala}{\mbox{$\cal{A}$}}

\newcommand{\calk}{\mbox{$\cal{K}$}}

\newcommand{\cals}{\mbox{$\cal{S}$}}

\newcommand{\reals}{\mbox{${\bf R}$}}
\newcommand{\rml}{\mbox{${\rm L}$}}

\newcommand{\NN} {\mbox{${\bf N}$}}

\newcommand{\cupnf}{\cup_{n=1}^{\infty}}

\newcommand{\sumnf}{\sum_{n=1}^{\infty}}
\newcommand{\sumin}{\sum_{i=1}^{n}}

\newcounter{excounter}

\begin{document}           
\newcommand{\sprt}{\vspace{0.2in}}
\newcommand{\sprth}{\vspace{0.1in}}
\newcommand{\beq}{\begin{eqnarray*}}
\newcommand{\eeq}{\end{eqnarray*}}
\newcommand{\iffp}{it follows from Proposition\ }
\newcommand{\Iffp}{It follows from Proposition\ }
\newcommand{\hardexercise}{
\hspace{-2.5mm}$ ^{*}$\hspace{1.5mm}}
\newcommand{\sbs}{\subset}
\newcommand{\sps}{\supset}
\newcommand{\stm}{\setminus}
\newcommand{\proof}{
  \sprt \bf Proof.  \rm \  }
\newcommand{\prooft}{
  \sprt \bf Proof of Theorems 9.17 and 9.18. \noindent \rm \  }
\newcommand{\thm}[1]{\stepcounter{props}
  \sprt\bf \noindent \thesection .\theprops .\
  #1:\rm \  }
\newcommand{\sol}[1]{
  \sprt\hspace{1cm}\bf  #1\par\hspace{-5pt}\nopagebreak\rm \  }
\newcommand{\exercises}{\sprt
\[ {\rm EXERCISES}\]\sprt\small
\setcounter{excounter}{0}\begin{list}%
{\thesection --\theexcounter}{\usecounter{excounter}
}}
\newcommand{\abst}{\sprt
\sprt
\begin{description}
\rightmargin=20 mm
}
\newcommand{\stareq}[3]{\vspace{7pt}$(#1)\hspace{#2cm}{\displaystyle
#3}$ \vspace{7pt}\\}
\newcommand{\bb}{
                 \item [ }
\large
\[ {\bf WEAKLY\; LINDELOF\; DETERMINED\; BANACH}\]
\[{\bf  SPACES\; NOT\; CONTAINING\;} \ell ^{1}(\NN )\]

\normalsize
\[{\rm by}\]
\[ {\rm Spiros\; A.\; Argyros}\]
\[ {\rm (Athens,\; Greece)}\]

\sprt\sprt
\abst
\item [\ \ \ \ \ \ \ \ ] {\bf ABSTRACT:}
{\em The class of countably intersected
families of sets is defined. For any such family we define
a Banach space not containing $\ell^{1}(\NN )$. Thus we
obtain counterexamples to certain questions related to
the heredity problem for W.C.G. Banach spaces. Among them we
give a subspace of a W.C.G. Banach space not containing
$\ell^{1}(\NN )$ and not being itself a W.C.G. space.}
\end{description}

\sprt
\sprt
{\bf INTRODUCTION}
In the present paper we deal with Banach spaces not containing
isomorphically the space $\ell^{1}(\NN )$.
The motivation for this study was a problem, posed to us by
S. Merkourakis, related to the heredity problem for weakly
compactly generated (W.C.G.) Banach spaces. The problem in question
is the following: Is every W.C.G. Banach space $X$ not containing
$\ell^{1}(\NN )$ hereditarily W.C.G.? It is well known by a classical
example due to Rosenthal [R] that the heredity problem for
W.C.G. has negative answer. On the other hand, M. Fabian has shown
in [F] that if the conjugate $X^{*}$ satisfies the Radon-Nikodym
property (R.N.P.) then the space $X$ is hereditarily W.C.G.
Actually M. Fabian proved the stronger result that if $X$ is weakly
countably determined (W.C.D.) and $X^{*}$ satisfies RNP
then $X$ is W.C.G. space. Later M. Valdivia in [V] extended this
result to the class of weakly Lindel\"{o}f determined (W.L.D.)
Banach spaces. (For the definitions we use here and related
results we refer the reader to [A-M-N], [A-M], [M-N] and
[D-G-Z]). For alternative proofs of these results we refer
the reader to [O-V-S] or [S] where the topological analogue of
the above results is also proved; that is, that every Corson
compact set which is also Radon-Nikodym set is Eberlein compact.

A class extending the W.L.D. Banach spaces with conjugate space
satisfying the RNP is W.L.D. spaces not containing
$\ell^{1}(\NN )$ and it is natural to ask for the solutions of
the same problems for this wider class. Our aim is to give
counterexamples to most of these questions. Thus in the
first section we introduce the {\bf countably intersected} families
of sets. This class of families includes the family of segments
of the dyadic tree and more generally the segments of any tree
whose each branch is a countable set. The most interesting
example of such a family is Recni\v{c}enko's space, which is
presented in detail in the first section of the paper.

In the second section we introduce the James norm for a family of
sets, which is the natural extension of James' norm for
James' tree space [J]. We prove that if the family is countably
intersected, then the resulting Banach space $X$ does not
contain $\ell^{1}(\NN )$. From this we conclude that \txs
a weakly \calk -analytic space $X$ not containing
$\ell^{1}(\NN )$ \st $X$ is not a subspace of a W.C.G. Banach
space as well as a W.L.D. Banach space $X$ not cobtaining
$\ell^{1}(\NN )$ and not being W.C.D. space. Thus we cannot
extend the results of Fabian and Valdivia to the class of
Banach spaces that do not contain $\ell^{1}(\NN )$.

The third section of the paper contains certain results related
to subspaces of W.C.G. spaces. By Rosenthal's example
\txs a probability measure $\mu$ \st the space $\rml^{1}(\mu )$
has a subspace which is not W.C.G.. Here we introduce the
{\bf quasi-Eberlein} sets and with each one of them we
associate a Banach space $X$ which is a subspace of a W.C.G.
space but it is not itself a W.C.G. space. We thus produce
counterexamples  to the heredity problem for W.C.G. Banach spaces
which are different from Rosenthal's space. It should be noted
that Merkourakis has recently shown that under ${\rm MA}$
every weakly \calk -analytic Banach space $X$ with
${\rm dim}X<2^{\omega}$ is a W.C.G. space. Thus it is not
possible to find a Z.F.C. counterxample to the heredity problem
for W.C.G. spaces of dimension $\omega_{1}$. Finally we give a Banach
space $X$ not containing $\ell^{1}(\NN )$ which is a subspace of
a W.C.G. Banach space but it is not a W.C.G. space.

\sprt
I wish to thank S. Merkourakis, A. Tsapalias and Th. Zachariades
for the useful discussions I had with  during the preparation
of the paper.

\sprt\sprt

\section{Countably intersected families}

This section is devoted to the definition and certain examples
of countably intersected families.

\thm{Definition}
A family \cals\ of subsets of a set $\Gamma $ is said to be
{\bf countably intersected} (C.I.) if the following conditions
are satisfied.

\newcounter{lista}
\begin{list}{(\alph{lista})  }{\usecounter{lista}}
\item \cals\ is a pointwise closed family of countable subsets
of $\Gamma $ containing all singletons.
\item For every $s,t$ in \cals , $s\stm t$ is the finite union
of \pair elements of \cals .
\item For every $t$ in \cals\ \txs $s_{t}$ in \cals\ containing
$t$ and such that \fe $s, t_{1},\ldots ,
t_{n}$ in \cals\ \txs $t$ in \cals , $t\sbs s\stm
\cup_{i=1}^{n} s_{t_{i}}$ and
\[ s\stm\left(\bigcup_{i=1}^{n} s_{t_{i}}\cup t\right)
\;\;{\rm is \; a\; finite\; set.}\]
\item For every $s$ in \cals\ the set
$L_{s}=\{ s\cap t:t\in\cals\} $ is countable.
\end{list}

\thm{Remarks} 
(i) Condition (a) in Definition 1.1 refers to the pointwise
closeness of the family $\{\chi_{s}:s\in\cals\}$ as a
subset of $\{ 0,1\}^{\Gamma }$. Since every $s$ in \cals\ is a
countable set, \cals\ is a Corson compact set ([A-M-N]).

(ii) Easy examples of C.I. families are pointwise closed families
containing only finite sets. These are called strong Eberlein
compact. We are interested in C.I. families lying beyond
the class of Eberlein compact sets. Next we will present certain
such examples.

\thm{Lemma} 
If $\{ s_{i}\}_{i=1}^{n}$ is a subfamily of a C.I. family \cals ,
\txs a \pair subfamily $\{ t_{j}\}_{j=1}^{l}$ of \cals\ such that
each $t_{j}$ is contained in some $s_{i}$ and
$\cup_{i=1}^{n}s_{i}=\cup_{j=1}^{l}t_{j}$.

\proof We prove it by induction on the length $n$ of the
family $\{ s_{i}\}_{i=1}^{n}.$

For $n=1$ it is clear. Suppose that $\cup_{i=1}^{n}s_{i}=
\cup_{j=1}^{l}t_{j}$ and $\{ t_{j}\}_{j=1}^{l}$  satisfies the
properties listed in the statement of the lemma. Then for
a given $s_{n+1}$,
\[ \bigcup_{i=1}^{n+1}s_{i}=\left(\bigcup_{j=1}^{l}t_{j}\right)\cup s_{n+1}
=\bigcup_{j=1}^{l}t_{j}\cup (\cdots ((s_{n+1}\stm t_{1})\stm t_{2})
\cdots\stm t_{l}).\]
Since \cals\ is a C.I. family, condition (b) in Definition
1.1 allows us to replace $s_{n+1}\stm t_{1}$ by a disjoint
union of elements of \cals , say $\cup_{j=1}^{l_{1}}t^{1}_{j}$.
Then $(s_{n+1}\stm t_{1})\stm t_{2}=\cup_{j=1}^{l_{1}}
(t_{j}^{1}\stm t_{2}) $ and again we can replace $t_{j}^{1}\stm t_{2}$
by a disjoint union of elements of \cals . Continuing in
this manner, we write $s_{n+1}\stm\cup_{j=1}^{l}t_{j}$ as a
disjoint union of elements of \cals .\hfill Q.E.D.

\[ {\bf Examples\;\; of\;\; countably\;\; intersected\;\; families}\]
\indent Since all the examples we will give below are related to
trees, we briefly recall some basic definitions and notation
 for trees.

A {\bf tree} is a partially ordered set $(T,\leq )$
\st for $\alpha\in T$ the set $\{ \beta \in T:\beta <\alpha\}$ is
well ordered.
The {\bf segments} of $T$ are defined as follows;
for $\beta \leq \alpha$ the segment $[\beta ,\alpha ]$
consists of all $\gamma \in T$ \st $\beta \leq \gamma \leq \alpha$.
The segment $[\beta ,\alpha )$ is defined similarly. If
$T$ has a unique minimal element $\alpha_{0}$, an {\bf initial segment}
is a segment of the form $[\alpha_{0},\alpha ]$ with $\alpha\in T$.
A {\bf branch} is any maximal segment. Finally, {\bf tail}
is any segment of the form $b\stm s$, where $b$ is a branch
and $s$ an initial segment contained in $b$.

\sprth
{\bf A.\ \  The C.I. family for the dyadic tree}

\sprth
The dyadic tree is the set
\[ \Delta =\bigcup_{n=0}^{\infty}\left\{ (n,k):
k=0,1,\ldots ,2^{n}-1\right\}\]
ordered by the relation $(n,k)<(n+1,l)$ if and only if
$l=2k$ or $l=2k+1$.

It is easy to see that the family \cals\ of all branches
and all tails of $\Delta $ defines a countably intersected family.
To check condition (c) of Definition 1.1 one should use as
$s_{t}$ the set $t$.

\sprth
{\bf B.\ \  The C.I. family for Todorcevi\v{c} tree}

\sprth
Consider a stationary subset $A$ of $\omega _{1}$ ($\omega _{1}$ is the
first uncountable ordinal number) \st $\omega _{1}\stm A$ is
also stationary. We set $T$ the set of closed subsets
of $A$. Recall that a  subset $K$ of $A$ is closed if and only if it
is a closed subset of $\omega _{1}$ with the topology induced
by the order. We define a partial order on $T$ by the rule:
$K_{1}\leq K_{2}$ if and only if $K_{1}$ is an initial segment of $K_{2}$.
Since, by definition, every stationary subset of $\omega _{1}$
meets every closed uncountable subset of $\omega _{1}$ and $A$,
$\omega _{1}\stm A$ are both stationary, we have that every branch of $T$ is
countable. The countably intersected family for
Todorcevi\v{c} tree is the family \cals\ of all segments of $T$.
Again, to check that this family is a C.I. family is
straightforward, the only point we want to note is that
for $t$ in \cals , the $s_{t}$ appearing in condition (c) of
Definition 1.1 is any initial segment containing $t$.

{\bf Remark.} As we defined the C.I. family for
Todorcevi\v{c} tree we could define a C.I. family for any
tree $T$ every branch of which is countable.

We turn now to Recni\v{c}enko's space which is more
complicated than the previous examples. We proceed to
a detailed presentation of it since, to the best of our
knowledge, there is no English reference for this example.
A brief presentation is given in [Ny].

\sprt
{\bf B.\ \  The C.I. family for Recni\v{c}enko's space}

\sprth
We set $\Gamma =\omega _{1}\times 2^{\omega }$, where $\omega _{1}$ and
$2^{\omega }$ are considered as ordinal numbers. There exists
a \sq $(\cala _{n})_{n\in {\bf N} }$, where each $\cala _{n}$
is a subset of $\Gamma $, satisfying the following properties.

\begin{list}{(\alph{lista})  }{\usecounter{lista}}
\item $\cala _{n}$ is a tree of height $\omega $ and, if we denote
by $\cala _{n}^{k}$ the (incomparable) elements of $\cala _{n}$
that belong to level $k$ of the tree, the following
properties are fulfilled.

\subitem (i) $\cala^{0}_{n}=\{ (0,n)\}$.

\subitem (ii) If $(\xi ,t)\in\cala^{k}_{n}$ then the set $(\xi,t)^{+}$
of the immediate successors of $(\xi,t)$ has cardinality
$2^{\omega }$ and every $(\zeta ,t)$ in $(\xi,t)^{+}$ satisfies
the condition $\xi <\zeta$.

\subitem (iii) For every $(\xi ,t)$ in $\cala _{n}$, $\xi >0$, \txs
exactly one $(\zeta ,s)$ in $\cala _{n}$ with
$(\xi ,t)\in (\zeta ,s)^{+}$.
\item The trees $\cala _{n}$, $n<\omega $, are connected by the
following preperties.

(i) If $A_{n}$ is an initial segment of $\cala _{n}$,
$A_{m}$ is an initial segment of $\cala _{m}$ and $n\neq m$
then $\# A_{n}\cap A_{m}\leq 1$.

(ii) For every $I$ infinite subset  of \NN ,
$A_{n}$ initial segments of $\cala _{n}$, $n\in I$, \st
$A_{n}\cap A_{m}=\emptyset$ \tx
uncountably many $(\xi ,t)$ in $\Gamma $ \st each
$(\xi ,t)$ extends $A_{n}$ for all $n$ in $I$.
\end{list}

\sprth
{\bf Construction of the trees $(\cala _{n})_{n\in {\bf N} }$}

\sprth
Each $\cala _{n}$ will be realized as $\cup_{\xi<\omega _{1}}
\cala^{\xi}_{n}$ with $ \cala^{\xi}_{n}$ satisfying the
following inductive properties.

\begin{list}{(\roman{lista})  }{\usecounter{lista}}
\item $\cala^{\xi}_{n}$ is a tree and $\cala^{\xi}_{n}
\sbs [0,\xi ]\times 2^{\omega }$.
\item For $\zeta<\xi<\omega _{1}$, $\cala^{\zeta}_{n}$ is a
subtree of $\cala^{\xi}_{n}$.
\item If $\xi =\zeta +1$, $I$ is an infinite subset of \NN\
and, for $n\in I$, $A_{n}$ is an initial segment of
$\cala^{\zeta}_{n}$ \st $A_{n}\cap A_{m}=\emptyset$ then
\txs $t<2^{\omega }$ \st $(\xi,t)$ extends $A_{n}$ for all
$n$ in $I$.
\item If $A_{i}$ is an initial segment of $A^{\xi}_{i}$,
$i=n,m$, $n\neq m$, then $\# A_{n}\cap A_{m}\leq 1$.
\end{list}

It is clear that the inductive properties listed above
imply that the family $\{ \cala _{n}\}_{n<\omega }$ satisfies
the desired properties.

\sprth
{\bf The inductive construction}
\sprth

We set $\cala^{0}_{n}=\{ (0,n)\}$.

If $\xi$ is a limit ordinal and $\cala^{\zeta}_{n}$
has been constructed for all $\zeta <\xi$, we set
$\cala^{\xi}_{n}=\cup_{\zeta <\xi}\cala^{\zeta}_{n}$.

If $\xi =\zeta +1$ and $\cala^{\zeta}_{n}$ has been
constructed, we consider all sequences $K=(A_{n})_{n\in I}$
\st $I$ is an infinite subset of \NN , $A_{n}$ is an
initial segment of $\cala _{n}$ and, for $n\neq m$,
$A_{n}\cap A_{m}=\emptyset$. The cardinality of all such $K$
is $2^{\omega }$ and we order them as
$\{ K^{\xi}_{t}\}_{t<2^{\omega }}$. W define $\cala^{\xi}_{n}$
as follows.
\[ \cala^{\xi}_{n}=\cala^{\zeta}_{n}\cup\left\{ (\eta,t)^{+}_{\xi}
:(\eta,t)\in\cala^{\zeta}_{n}\right\} ,\]
where $(\eta,t)^{+}_{\xi}=\{ (\xi,t'):{\rm there\; exists}\;
A\in K^{\xi}_{t'}, A\sbs\cala _{n}\;{\rm and}\;
\max A=(\eta ,t)\}$. We order $\cala^{\xi}_{n}$ in the
obvious way.

Thus the inductive construction is complete.

{\bf Reecni\v{c}enko's space} is the following family.
\[ R=\{ s:s{\rm\;is\; a\; segment\; for\; some\;} \cala _{n}\} .\]
It is easy to check that $R$ is a family of countable
subsets of $\Gamma $ which is countably intersected.

\sprth
\thm{Proposition}
{\em $R$ is a Talagrand compact set which is not Eberlein
compact.}

\sprth
The proof will follow from the next two lemmata.

Before we state the first lemma we need some notation.
We denote by $FS(\NN )$ the (countable) set of all finite
sequences of natural numbers. Also, for $d\in \NN^{{\bf N}}$,
we denote by $dn$ the first $n$ terms of the \sq $d$.

\thm{Lemma}
{\em There exists a partition $\{ \Gamma _{\phi }:\phi \in FS(\NN )\}$
of $\,\Gamma $ \st \fe $d\in{\bf N}^{{\bf N}}$ and every
$(\xi _{n},t_{n})\in \Gamma _{\phi n}$ the set $\{ (\xi _{n},t_{n}
)\}_{n\in{\bf N}}\cap s$ is finite for all $s$ in $R$.}

\sprt
It is well known that every pointwise closed family
satisfying the conclusion of the lemma defines a Talagrand
compact set (see [A-N]).

\proof
For $n\in{\bf N}$, $k=0,1,2,\ldots$, we define
$\Gamma ^{n,k}=\{ (\xi,t):(\xi,t)\in\cala^{k}_{n}\}$
(recall $\cala^{k}_{n}$ is the ${\rm k}^{th}$ level of $\cala _{n}$)
and for $\phi \in FS(\NN )$, $\phi =(n_{1},\ldots ,n_{s}),$
we set \[ \Gamma _{\phi }=\bigcup_{i=1}^{s}\Gamma ^{i,n_{i}}.\]
It is easy to check that $\{ \Gamma _{\phi }:\phi \in FS(\NN )\}$
satisfies the desired properties.

\vspace{-0.05in}
\hfill Q.E.D.

\sprt
The following lemma will show that $R$ is not an Eberlein compact set.
Indeed, if a totally disconnected set is an Eberlein compact set,
as S. Merkourakis noticed, the algebra of clopen subsets is a
$\sigma $-weakly compact set. This follows from Rosenthal's
characterization of the Eberlein compacta [R].
Therefore, if $R$ were an Eberlein compactum and we considered
$R$ as a subset of $\{ 0,1\}^{\Gamma }$, we could find a
partition of $\Gamma $ into a countable family $\Gamma _{n}$ \st \fe
$\{ (\xi _{k},t_{k})\}_{k\in{\bf N}}\sbs \Gamma _{n}$ \txs a sub\sq
$\{ \Gamma (\xi _{k},t_{k})\}_{k\in M}$ converging pointwise to zero.
 But this contradicts the conclusion of the next lemma. For
 alternative proofs see [So], [Fa].

\thm{Lemma}
{\em For every partition $(D_{i})_{i\in{\bf N}}$ of the set $\Gamma$ \txs
$M\in R$ and $i_{0}\in{\bf N}$ \st $M\cap D_{i_{0}}$ is
an infinite set.}

\proof We define the following set
\beq L\!\!&=\!\!&\{ i\in{\bf N} :{\rm there\; exists\; an\; infinite\;
countable\; family\;} \{M_{k}\}_{k\in{\bf N}}\; {\rm of\;
pairwise}\\
& &\;\;\;\; {\rm disjoint\; elements\; of\; \cala\; and\;
(\xi ,t)\in D_{i}\; such\; that}\; M_{k}\cup\{ (\xi ,t)\}\in
\cala\\
& &\;\;\;\; {\rm for\; all}\; k\in{\bf N} \}.
\eeq
It is clear that $L$ is non empty set. For each $i\in L$
and some $(M_{k})_{k=1}^{\infty}$ witnessing the belongingness of
$i$ to $L$ we set
\[ E_{i}=\{ n_{k}:M_{k}\in\cala_{n_{k}}\} .\]
Each $E_{i}$ is an infinite set.

Assuming the failure of the conclusion, we examine the following
two mutually exclusive cases.

{\em Case 1.} The set $L$ is an infinite set.

Then we could define $\phi :L\lra\NN$ one-to-one function \st
for all $i\in L$, $\phi (i)\in E_{i}$. Notice that our assumptions
imply that \fe $i\in L$ \txs $M_{i}\in\cala_{\phi (i)}$ \st
$M_{i}\cap D_{i}\neq\emptyset$, $\max M_{i}\in D_{i}$ and $M_{i}$ does
not have extensions $(\xi ,t)$ with $(\xi ,t)\in D_{i}$.
But the family $\{ M_{i}\}_{i\in L}$ consists of \pair
elements of \cala , hence it has at least one extension;
that is, \txs $(\xi ,t)\in \Gamma $ with $M_{i}\cup\{ (\xi,t)
\}\in\cala$ for all $i\in L$. Since $\{ D_{i}\}_{i\in{\bf N}}$
is a partition of $\Gamma $, \txs $i_{0}$ in \NN\ with
$(\xi,t)\in D_{i_{0}}$. Further $i_{0}\in L$, hence
$(\xi,t)$ extends $M_{i_{0}}$ which contradicts the maximality
of $M_{i_{0}}$ and the proof in this case is complete.

{\em Case 2.} The set $L$ is a finite set.

In this case we produce inductively a set $M$ \st $M\in R$
and $M\cap D_{i}$ is an infinite set for some $i\in L$.

Indeed, we start with the family $\{ ((0,n))\}_{n\in{\bf N}}$.
Then \txs $i_{1}\in L$ and an infinite set $\{ (\xi ^{1}_{n},t^{1}_{n})
\}_{n\in{\bf N}}\sbs D_{i_{1}}$ \st
$\left( (0,n), (\xi ^{1}_{n},t^{1}_{n})\right)
\in\cala$ and $(\xi ^{1}_{n},t^{1}_{n})\neq (\xi ^{1}_{m},t^{1}_{m})$ for
$n\neq m$. In the next step we consider the family
$\{ ((0,n),(\xi ^{1}_{n},t^{1}_{n}))\}_{n\in{\bf N}}$ and we select
$i_{2}\in L$, $\{(\xi ^{2}_{n},t^{2}_{n})\}_{n\in{\bf N}}\sbs D_{i_{2}}$
\st \[ \left( (0,n),(\xi ^{1}_{n},t^{1}_{n}),(\xi ^{2}_{n},t^{2}_{n})\right)
\in\cala\;\;{\rm and}\;\; (\xi ^{2}_{n},t^{2}_{n})\neq
(\xi ^{2}_{m},t^{2}_{m})\; {\rm for}\; n\neq m.\]
We proceed inductively in the same manner. Since $L$ is a finite set,
we find $i_{0}\in L$ and $(n_{k})_{k\in{\bf N}}$ an infinite
set \st for all $n\in{\bf N}$, $((\xi^{n_{k}}_{k},
t^{n_{k}}_{k})\in D_{i_{0}}$. Then clearly for all $n\in{\bf N}$,
$\left( (0,n),(\xi ^{1}_{n},t^{1}_{n}),\ldots ,(\xi ^{k}_{n},t^{k}_{n}),
\ldots \right)\cap D_{i_{0}}$ is an infinite set, a contradiction
and the proof of the lemma is complete.\hfill Q.E.D.

\sprth
Actually the space $R$ satisfies the following stronger condition
that we will use in the third section.

\thm{Lemma}
{\em For every $(\Gamma _{d})_{d\in D}$ partition of $\,\Gamma $  into
\pair countable sets and every $(D_{i})_{i\in {\bf N}}$ countable
partition of $\Gamma $ \txs $M\in R$ and $i_{0}\in {\bf N}$ \st
$M\cap D_{i_{0}}$ is an infinite set and $\# M\cap \Gamma _{d}\leq 1$
for all $\delta\in D$.}

\proof The proof is essentially the same as in the
previous lemma. We only have to take into account the additional
property that $\# M\cap \Gamma _{d}\leq 1$. Thus we define the set
$L$ as:
\beq L\!\!&=\!\!&\{ i\in{\bf N} :{\rm there\; exists\; an\; infinite\;
countable\; family\;} \{M_{k}\}_{k\in{\bf N}}\; {\rm of\;
pairwise}\\
& &\;\;\;\; {\rm disjoint\; elements\; of\; \cala\; and\;
(\xi ,t)\in D_{i}\; such\; that}\; M_{k}\cup\{ (\xi ,t)\}\in
\cala\\
& &\;\;\;\; {\rm and}\; \#\, M_{k}\cup\{ (\xi ,t)\}\cap \Gamma _{d}
\leq 1\; {\rm for\; all}\; k\in{\bf N} \; {\rm and}\; d\in D \}.
\eeq
Using conditions (b)-(ii) of the definition of the families
$(\cala _{n})_{n\in {\bf N}}$, we proceed with the proof in the same
manner as in Lemma 1.6. \hfill Q.E.D.

\section{Banach spaces not containing $\ell ^{1}(\NN)$}

Let $\Gamma $ be a non empty set and \cals\ a family of subsets
of $\Gamma $. We denote by $\Phi $ the linear space of all finitely
supported real valued functions defined on $\Gamma $. For $a\in
\Gamma $ we denote by $e_{a}$ the characteristic function of the
one point set $\{ a\}$. Clearly the family
$\{ e_{a}\}_{a\in \Gamma }$ is a Hamel basis for the space $\Phi $.

\thm{Definition}
Let $\Gamma $ be a non empty set and \cals\ a family of subsets of
$\Gamma $. For $\phi $ a finitely supported real valued function on
$\Gamma $, we define the {\bf James-\cals} or $J-\cals $ norm as:
\[ \ \phi \ =\sup\left(\sumin\left(\sum_{a\in s_{i}}
\phi (a)\right) ^{2}\right)^{1/2},\]
where the supremum is taken over all families
$\{ s_{1},s_{2},\ldots ,s_{n}\}$ of \pair elements of \cals .  We
denote by $J\cals$ the completion of $\Phi $ in the above defined norm.

Our goal is to prove the next result.

\thm{Theorem} {\em For every C.I. family \cals\ the space $J\cals$
does not contain isomorphically $\ell ^{1}(\NN )$.}

\sprt
{\bf Remark.} (i) This result is known for certain C.I. families.
For example, the classical James Tree space [J] is the space
defined by $J-\cals$ norm for the C.I. family of the dyadic tree.
(More precisely, the family used in the original James'
definition is that \cals\ that consists of all segments of
the dyadic tree.

(ii) The proof of the theorem is divided into two parts. In the
first we will give a representation of the extreme points of the
unit ball of the dual $J\cals ^{*}$. In the second part we show
that every bounded \sq in $J\cals$ has a weak Cauchy sub\sq.
For this we use the fact that the family \cals\ is countably
intersected together with the representation of the extreme
points of the first step and Rainwater's theorem.

\sprth
 We begin with some notation and definitions.

With each $s$ in \cals\ we associate a functional $s^{*}$ in
$J\cals ^{*}$ defined by
\linebreak
$s^{*}(\phi )=\sum_{\alpha\in s}\phi (\alpha )
\;\;$
\fe $\phi $ finitely supported real valued function. It follows
from the definition of $J\cals$ norm that $s^{*}(\phi )\leq
\ \phi $, hence $s^{*}\ \leq 1$ and if $s\neq\emptyset$ then
$s^{*}\ =1$. We also define the set
\[ D=\left\{\sumin \lambda _{i}s_{i}^{*}:\lambda _{i}\in\reals, s_{i}\in\cals,
\{ s_{i}\}_{i=1}^{n}\;{\rm are\; pairwise\; disjoint\; and}\;
\sumin \lambda ^{2}_{i}\leq 1\right\} .\]
A simple argument shows that for $\{ s_{i}\}_{i=1}^{n}$ pairwise
disjoint \[ \left( \sumin \lambda _{i}s_{i}^{*}\right)\leq
\left(\sumin \lambda ^{2}_{i}\right)^{1/2},\]
hence $D$ is a subset of the unit ball of $J\cals ^{*}$.

\thm{Lemma}
{\em The $w^{*}$ closure of $D$ contains the extreme points of the
unit ball of $J\cals ^{*}$.}

\proof Suppose that the conclusion is false. Then \txs an extreme
point $x^{*}$ of $B_{J{\cal S} ^{*}}$ with $x^{*}\not\in
\overline{D}^{w^{*}}$. Hence we could find a $w^{*}$-neighborhood of
$x^{*}$ in $B_{J{\cal S} ^{*}}$ disjoint from $\overline{D}^{w^{*}}$.
But since the $w^{*}$ slices of $x^{*}$ define a neighborhoods basis
for the $w^{*}$ topology, \txs a slice $S(x^{*},x,t)$ disjoint from $D$.
Further we may assume that $x$ is a finitely supported function
on $\Gamma $ with $\ x\ =1$. From all the above we get that

(i) For some $\epsilon >0$, $x^{*}(x)>\sup_{y^{*}\in D}y^{*}(x)+\epsilon $.

(ii) There exists $\{ s_{i}\}_{i=1}^{n}$ \pair and \st
\[ 1=\ x\ =\left(\sumin\left(\sum_{a\in s_{i}}x(a)
\right)^{2}\right)^{1/2}.\]
We set $\lambda _{i}=\sum_{a\in s_{i}}\phi (a)$ and
$y^{*}=\sumin \lambda _{i}s_{i}^{*}$.
Then $y^{*}\in D$ and
\[ 1=y^{*}(x)<x^{*}(x)-\epsilon \leq 1-\epsilon ,\]
a contradiction and the proof is complete.\hfill Q.E.D.

\sprt
In the sequel we identify each $s\in\cals$ with its characteristic
function on $\{ 0,1\}^{\Gamma }$, and \cals\ with the corresponding
(closed) subset of $\{ 0,1\}^{\Gamma }$. Recall that for $a\in \Gamma $,
$e_{a}$ denotes $\chi_{\{ a\} }$ and note that $\ e_{a}\ =1$
and the family $\{ e_{a}\}_{a\in \Gamma }$ generates the space
$J\cals$. Finally, we denote by $B_{\ell ^{2}}$ the unit
ball of $\ell ^{2}(\NN )$ endowed with the weak topology.
We set
\[ V=\left\{ (\lambda _{n})_{n\in{\bf N}}\sbs B_{\ell ^{2}}:
( \lambda _{n})_{n\in{\bf N}}\;{\rm  is\; a\; decreasing\;
 sequence}\,\right\} ,\]
which is a closed subset of $B_{\ell ^{2}}$.

\thm{Lemma}
{\em Let $E\sbs S^{{\bf N}}$ be defined by $(s_{n})_{n\in{\bf N}}\in E$
 if and only if $(s_{n})_{n\in{\bf N}}$ is a \pair family. Then

(i) $E$ is closed subset of $S^{{\bf N}}$.

(ii) The function
\[ T:V\times E\lra\left(\overline{D}^{w^{*}},w^{*}\right)\;
de\!f\!ined\; by\; T\left( (\lambda_{n})_{n\in{\bf N}},
(s_{n})_{n\in{\bf N}}\right) =\sumnf\lambda_{n}s_{n}^{*}\]
is continuous and onto.}

\proof
(i) It is straightforward and it is left to the reader.

(ii) Let $x_{i}=\left( (\lambda_{n,i})_{n\in{\bf N}},
(s_{n,i})_{n\in{\bf N}}\right)$ and
$x=\left( (\lambda_{n})_{n\in{\bf N}},
(s_{n})_{n\in{\bf N}}\right)$ and suppose that the net
$(x_{i})_{i\in I} $ converges to $x$.

We will show that $Tx_{i}$ converges in the $w^{*}$ topology
to $Tx$. For this it is enough to show that
\[ Tx_{i}(e_{a})\lra Tx(e_{a})\;\;{\rm for
\; all\;} a\in \Gamma .\]
Indeed, there are two cases.

{\em Case 1.} $a\in\cup_{n\in{\bf N}}s_{n}$.
Then there exists a unique $n_{0}\in{\bf N}$ \st $a\in s_{n_{0}}$.

Since $E\sbs \left(\{ 0,1\}^{\Gamma }\right)^{{\bf N}}$, denoting
by $(a,n_{0})$ the $a^{{\rm th}}$ coordinate of $\Gamma $
appearing in the $n_{0}$ copy of it, we have
$x_{i}(a,n_{0})\lra x(a,n_{0})$, hence \txs $i_{0}$ \st
for all $i\geq i_{0}$, $x_{i}(a,n_{0})= x(a,n_{0})=1$ and
$x_{i}(a,m)=0$ for all $m\in{\bf N}$, $m\neq n_{0}$.
Consequently, for $i\geq i_{0}$, $Tx_{i}(e_{a})=
\lambda _{(n_{0},i)}$ which converges to $\lambda _{n_{0}}=
Tx(e_{a})$.

{\em Case 2.} $a\not\in \cup_{n\in{\bf N}}s_{n}$.

In this case we will show that $Tx_{i}(e_{a})$ converges
in the $w^{*}$ topology to 0.

Indeed, let $\epsilon >0$. Since \fa $i\in I$, $\sumnf
\lambda ^{2}_{(n,i)}\leq 1$ and $(\lambda _{(n,i)})_{n\in{\bf N}}$ is
decreasing, \txs $n_{0}\in{\bf N}$ \st for all $i\in I$ and
$m>n_{0}$, $\lambda _{(m,i)}<\epsilon $.

Also, since $a\not\in\cup_{n\in{\bf N}}s_{n}$, \txs $i_{0}$
\st $a\not\in s_{(k,i)}$ for all $i\geq i_{0}$ and
$k=1,2,\ldots ,n_{0}$.

Then for $i\geq i_{0}$, $Tx_{i}(e_{a})\leq\sup
\{ \lambda _{(m,i)}:m>n_{0}\}\leq \epsilon $ and this proves the result
in the second case.\hfill Q.E.D.

\sprth
We turn to some consequences of the previous results.

\thm{Corollary}
{\em Every extreme point $x^{*}$ of the unit ball of
$J\cals ^{*}$ is represented as
\[ x^{*}=\sumnf \lambda _{n}s_{n}^{*},\]
where $\sumnf \lambda ^{2}_{n}=1$, $(s_{n})_{n\in{\bf N}}$ is a \pair
\sq of elements of \cals\ and the series converges in the
norm topology of $J\cals ^{*}$.}

\thm{Corollary}
{\em If the family \cals\ defines a Talagrand compact set, then
the space $J\cals$ is  weakly $\calk$-analytic.}

\proof
It is well known that the class of Talagrand compact sets is
closed \wrt closed subsets, countable products and continuous
images [T]. Hence $\overline{D}^{w^{*}}$ is a Talagrand
compact set containing the extreme points of $B_{J{\cal S} ^{*}}$.
Therefore $J\cals$ is isometric to a closed  subspace of
the space $C(\overline{D}^{w^{*}})$ and the
proof is complete.\hfill Q.E.D.

\thm{Corollary}
Let $(x_{n})_{n\in{\bf N}}$ be a bounded \sq in $J\cals$. If
\fa $s$ in \cals\ the \sq $(s^{*}(x_{n}))_{n\in{\bf N}}$ is
Cauchy, then $(x_{n})_{n\in{\bf N}}$ is weakly Cauchy.

\proof
It follows from Corollary 2.5 that \fe $x^{*}$ extreme point of
$B_{J{\cal S} ^{*}}$, $(x^{*}(x_{n}))_{n\in{\bf N}}$ is Cauchy. The result
now follows from Rainwater's theorem [Rai].\hfill Q.E.D.

\sprt
The second part of the proof of the main theorem is contained
in the next lemma.

\thm{Lemma}
{\em Every norm bounded \sq $(\phi _{n})_{n\in{\bf N}}$ in $J\cals$
consisting of finitely supported functions contains a sub\sq
$(\phi _{n})_{n\in B}$ \st for all $s\in\cals$,
$(s^{*}(\phi _{n}))_{n\in B}$ is Cauchy.}

\proof
Assume that \fa $n\in{\bf N}$, $\ \phi _{n}\\leq 1$. We claim
the following.

{\em Claim 1.} For every $\epsilon >0$ and every $(\phi _{n})_{n\in B}$
sub\sq of the given \sq \txs a finite set
$\{ s_{1},\ldots ,s_{k}\}$ subset of \cals\ and an infinite $B'$
subset of $B$ \st \fe $s\in\cals$ with $s\cap s_{i}=\emptyset$,
$i=1,\ldots ,k$, we have
\[ \limsup_{n\in B'}s^{*}(\phi _{n})\leq \epsilon .\]
To see this, we begin with an infinite $B$ subset of \NN\
and we choose $s_{1}$ in \cals , if such exists, so that
\[ \limsup_{n\in B}s_{1}^{*}(\phi )> \epsilon .\]
Then we choose $B_{1}$ infinite \st
\[ s^{*}_{1}(\phi _{n})>\epsilon \;\;{\rm for\; all}\; n\in B_{1}.\]
If the pair $\{ s_{1}\} ,B_{1}$ does not satisfy the conclusion
of the claim, we repeat the procedure to find an $s_{2}$ in
\cals\ disjoint from $s_{1}$ \st
\[ \limsup_{n\in B_{1}}s_{2}^{*}(\phi _{n})> \epsilon .\]
Choose $B_{2}$ subset of $B_{1}$ \st
\[ s^{*}_{2}(\phi _{n})>\epsilon \;\;{\rm for\; all}\; n\in B_{2}.\]
Notice that \fa $n\in B_{2}$, $\ \phi \>\epsilon \sqrt{2}.$
Therefore, if $k=\min\{ n:\epsilon \sqrt{n}>1\}$ and we
repeat the procedure at most $k$ times, we will get
$s_{1},s_{2},\ldots ,s_{l}$, for some $l\leq k$ and an
infinite $B'$ subset of $B$ \st the conclusion of the
claim is satisfied.

Next we apply the claim infinitely many times to choose
a \dsq $(B_{m})_{m\in{\bf N}}$ of subsets of \NN, and families
$\{ s_{1,m},\ldots ,s_{k_{m},m}\}$ subsets of \cals\
\st the pair $(\phi _{n})_{n\in B_{m}},(s_{1,m},\ldots ,s_{k_{m},m})$
satisfies the conclusion of the claim for $\epsilon =1/m$.
Choose a \sq $(\phi _{n})_{n\in B'}$ which is almost
contained in every $B_{m}$.
The family \cals\ is countably intersected. Therefore,
for every $s_{i,m}$ \txs $\tilde{s}_{i,m}$ containing
$s_{i,m}$ and satisfying condition (c) of Definition 1.1.
Notice also that the set
\[ L=\bigcup_{m=1}^{\infty}\bigcup_{i=1}^{k_{m}}L_{\tilde{s}_{i,m}}\]
is countable (condition (d) of Definition 1.1).
Further, for $s,t$ in \cals , condition (b) of Definition 1.1
ensures that $(s\stm t)^{*}$ belongs to $J\cals ^{*}$, hence
$(s\cap t)^{*}=s^{*}\stm (s\stm t)^{*}$ is also in $J\cals ^{*}$.
Therefore, using a diagonal argument, we may choose a
sub\sq $(\phi _{n})_{n\in B''}$ of $(\phi _{n})_{n\in B'}$
\st $(t^{*}(\phi _{n}))$ converges for all $t^{*}\in L$ and
all $t^{*}$ finite subsets of $\cup_{n\in B'}{\rm supp} \phi _{n}$.

Notice that the last property implies that
$(f^{*}(\phi _{n}))_{n\in B''}$ converges for all $f$ in the
class of finite sets.

{\em Claim 2.} The \sq $(s^{*}(\phi _{n}))_{n\in B''}$
is Cauchy for all $s\in\cals$.

Indeed, consider any $s$ in \cals\ and $m\in{\bf N}$. We set
\[ s_{m}=s\stm\bigcup_{i=1}^{k_{m}}\tilde{s}_{i,m}.\]
Also, we choose $\{ t_{j,m}\}_{j=1}^{l_{m}}$ \pair
subfamily of \cals\ \st
$ \cup_{i=1}^{k_{m}}\tilde{s}_{i,m}\! =
\cup_{j=1}^{l_{m}}t_{j,m}$ (Lemma 1.3) and each
$t_{j,m}$ is contained in some $\tilde{s}_{i,m}$.
Finally, we set $d_{j,m}=t_{j,m}\cap s$ and we choose
$t_{m}\sbs s_{m}$ with $s_{m}\stm t_{m}$ being a finite set.
Then \beq \limsup_{n\in B''}s^{*}(\phi _{n})&-&
\liminf_{n\in B''}s^{*}(\phi _{n})=
\limsup_{n\in B''}t^{*}_{m}(\phi _{n})\\
&+&\lim_{n\in B''}\sum_{j=1}^{l_{m}}d^{*}_{j,m}(\phi _{n})+
\lim_{n\in B''}(s_{m}\stm t_{m})^{*}(\phi _{n})\\
&-&\left(\liminf_{n\in B''}t^{*}_{m}(\phi _{n})+
\lim_{n\in B''}d^{*}_{j,m}(\phi _{n})+
\lim_{n\in B''}(s_{m}\stm t_{m})^{*}(\phi _{n})\right)\\
&=&\limsup_{n\in B''}t^{*}_{m}(\phi _{n})-\liminf_{n\in B''}
t^{*}_{m}(\phi _{n})\leq\frac{1}{m}.
\eeq
Since the last inequality holds \fe $m\in{\bf N}$, we get
that $(s^{*}(\phi _{n}))_{n\in B''}$ is Cauchy.\hfill Q.E.D.

\sprth
{\bf Proof of Theorem 2.2.}
It is well known that a Banach space $X$ does not contain
isomorphically $\ell ^{1}(\NN)$ if every bounded \sq
$(x_{n})_{n\in{\bf N}}$ has a weakly Cauchy subsequence.
Given a bounded \sq $(x_{n})_{n\in{\bf N}}$ in $J\cals$, we
approximate each $x_{n}$ by some $\phi _{n}$, a function of finite
support \st $\ x_{n}-\phi _{n}\\leq\frac{1}{n}$. Applying
Lemma 2.8 and Corollary 2.7, we find a sub\sq
$(\phi _{n})_{n\in B}$ which is weakly Cauchy. Since
$x_{n}-\phi _{n}$ converges in norm to zero, we get that
$(x_{n})_{n\in{\bf B}}$ is also weakly Cauchy.\hfill Q.E.D.

\thm{Corollary}
{\em There exists a Banach space $X$ not containing
$\ell ^{1}(\NN )$ \st $X$ is weakly $\calk$-analytic
but not a subspace of a W.C.G.}

\proof
This is the $J\cals$ space for \cals\ the C.I. family for
Recni\v{c}enko's space. Indeed, in Corollary 2.6 we proved that
$J\cals$ is weakly $\calk$-analytic and the unit ball of
$J\cals^{*}$ contains a subset which is homeomorphic in the
$w^{*}$ topology to \cals . Hence it is not Eberlein compact
set and $X$ is not a subspace of a W.C.G.

\vspace{-0.05in}
\hfill Q.E.D.

\thm{Corollary}
{\em There exists a  Banach space $X$ not containing $\ell ^{1}(\NN )$
\st $(B_{X^{*}},w^{*})$ is Corson compact set, or equivalently
$X$ is a W.L.D. space,  and $X$ is not weakly
$\calk$-analytic.}

\proof
This is the $J\cals$ space for \cals\ the C.I. family for
Todorcevi\v{c} tree.
Indeed, \cals\ is not a Talagrand compact set and it is
homeomorphic to a closed subset of the unit ball of $J\cals^{*}$.
This shows that $J\cals$ is not weakly \calk -analytic. On the
other hand \cals\ has property (M) [A-M (Proposition 3.10)] and
hence $B_{J{\cal S}^{*}}$ in the $w^{*}$- topology is a Corson
compact set. \hfill Q.E.D.

\section{Subspaces of W.C.G. spaces}

In this section we give examples of subspaces of W.C.G. spaces
which are not W.C.G. The first part contains a general method
of producing such examples. In the second part we construct
an example of this kind that satisfies the additional property
that $\ell ^{1}({\bf N})$ is not isomorphically embedded into the space.

 We begin with the definition of a class of compact spaces.

\thm{Definition} A pointwise closed family of countable subsets
 of a set $\Gamma $ is said to be {\bf quasi-Eberlein} if the
following conditions are satisfied

\begin{list}{(\roman{lista})  }{\usecounter{lista}}
\item The singletons of $\Gamma $ are contained in \cals ,
and each $s$ in \cals\ is an increasing union of
finite sets $t$ from \cals .
\item $\cals =\cup_{n=1}^{\infty}\cals _{n}$, where each $\cals _{n}$
is an Eberlein compact set and $\Gamma $ is a subset of $\cals _{1}$.
\item For all partitions $(\Gamma _{d})_{d\in D}$,
$(\Gamma  _{n})_{n=1}^{\infty}$ of
$\Gamma $ \st $\Gamma _{d}$ are countable and
\pair \txs an $s$ in \cals\ and $n_{0}\in {\bf N}$ \st
$s\cap \Gamma _{n_{0}}$ is an infinite set while
$\# s\cap \Gamma _{d}\leq 1$ for all $d$ in $D$.
\end{list}

\thm{Remark}
(a) The first part of condition (i)
and the second part of condition (ii) are
not so important. If \cals\ does not satisfy these two conditions,
we may consider the family $\cals '=\cals\cup \Gamma $ and
$\cals _{1} '=\cals_{1}\cup \Gamma $ and the family $\cals '$ is a
quasi-Eberlein set.

(b) The fact in condition (iii) that \fe partition
$(\Gamma  _{n})_{n\in{\bf N}}$ of $\Gamma $ \txs
 $s$ in \cals\ and $n_{0}\in {\bf N}$ \st
$s\cap \Gamma _{n_{0}}$ is an infinite set shows, as we explained
in the case of Recni\v{c}enko's space, that \cals\ is not an Eberlein
compact set. The second part of condition (iii) will be used
to prove that certain Banach spaces defined by means of
 quasi-Eberlein sets are not W.C.G.

\thm{Lemma}
{\em If \cals\ is a quasi-Eberlein set and $(\cals _{n})_{n\in{\bf N}}$
its partition  into Eberlein compact sets, then the set
\[ \cals E=\left\{ \frac{1}{n}\chi_{s}:s\in\cals _{n}\right\}\]
is  Eberlein compact.}

\proof It is easy to see that $\cals E$ is a closed subset
of $[0,1]^{\Gamma }$ and using Rosenthal's criterion [R] we get the
desired result.

\thm{Definition}
We call the set $\cals E$ the {\bf Eberleinization} of the
quasi-Eberlein set \cals .

\thm{Examples}
We give two examples of quasi-Eberlein sets. Notice that by
a result due to Sokolov [S], every quasi-Eberlein set is a
Talagrand compact set.

\sprth
{\bf A. \ Talagrand's space:}
Recall the definition of Talagrand's example that is given in [T].
We set $\Gamma ={\bf N}^{{\bf N}}$; a subset $s$ of $\Gamma $
 is said to be
{\bf admissible} if and only if \txs $n\in {\bf N}$ \st for
$\gamma _{1},\gamma _{2}\in s$, $\gamma _{1}\neq \gamma _{2}$,
we have $\gamma _{1}(1)=
\gamma _{2}(1),\gamma _{1}(2)=\gamma _{2}(2),\ldots ,\gamma _{1}(n-1)=
\gamma _{2}(n-1)$ and $\gamma _{1}(n)\neq \gamma _{2}(n)$.
We call the number $n$ the {\bf characteristic} of the
set $s$. We also agree that if $\# s\leq 1$ then every
$n\in {\bf N}$ is the characteristic of $s$.
The family \cals\ of all admissible sets is pointwise
closed and if $s$ is admissible, $t$ a subset of $s$, then
$t$ is also admissible. It is an application of Baire's category
theorem that if $(\Gamma _{n})_{n\in{\bf N}}$ is a
partition of $\Gamma $, then
\txs an infinite admissible $s$ contained in some $\Gamma _{n}$.
Therefore, if $\{ \Gamma _{d}\}_{d\in D}$ is a partition of
$\Gamma $ into countable \pair sets then \txs a partition
$\{ \Delta _{k}\}_{k\in {\bf N}}$ of $\Gamma $ \st
$\# \Delta _{k}\cap \Gamma _{d}\leq 1$ for all $k\in {\bf N}$,
$d\in \Delta $.
Hence for given $\{ \Gamma _{d}\}_{d\in D},
\{ \Gamma _{n}\} _{n\in{\bf N}}$ partitions of $\Gamma $ as in Definition
3.1 (iii), we consider as above the partition
$\{ \Delta _{k}\}_{k\in {\bf N}}$ corresponding to the family
$\{ \Gamma _{d}\}_{d\in D}$ and the joint partition
$\{ \Gamma _{(n,k)}\}_{(n,k)\in {\bf N}\times {\bf N}}$, where
$\Gamma _{(n,k)}=\Gamma _{n}\cap \Delta _{k}$. It is easy to see that
every infinite admissible $s$ contained in some
$\Gamma _{(n,k)}$ satisfies condition (iii) of Definition 3.1.
Finally we set
\[ \cals _{n}=\{ s\in\cals : {\rm the\; characteristic \;
of}\; s \;{\rm is}\; n\} .\]
Then $\cals _{n}$ is an Eberlein compact set and
$\cals =\cupnf\cals _{n}$. Therefore \cals\ is a
 quasi-Eberlein compact set.

\sprth
{\bf B.\ Recni\v{c}enko's space} : We recall from the first section
of the paper that Recni\v{c}enko's space is the set
\[ R=\{ s : {\rm there\;  exists}\; n\in {\bf N}\; {\rm with}
\; s\;{\rm a\; segment\; of\; the\; tree}\;\cala _{n}\} .\]
Also, the order type of each branch of $\cala _{n}$ is
$\omega $. Therefore each $\gamma $ in $\cala _{n}$ belongs to some
$\cala _{n}^{k}$, where $k\in {\bf N}$ denotes the ${\rm k}^{th}$
level of $\cala _{n}$. Hence the correspondence
$F:\Gamma \cap\cala _{n}\lra c_{0}(\Gamma )$ defined as
$\cala _{n}^{k}\ni \gamma \lra\frac{1}{k+1}e_{\gamma }$ is extended
to a homeomorphism from the segments of $\cala _{n}$ to a
subset of $(B_{c_{0}(\Gamma )},w)$. Thus the set
\[ R_{n}=\{ s:s\;{\rm is\; a\; segment\; of}\;\cala _{n}\}\]
is an Eberlein compact set. We set
\[ R'=R\cup \Gamma \;\;{\rm and}\;\; R_{1}'=R_{1}\cup \Gamma .\]
Then $R'$ is a quasi-Eberlein compact set (see Lemma 1.7).

\thm{Theorem}
{\em Let \cals\ be a quasi-Eberlein set and $\cals E$ its Eberleinization
that is naturally contained in the the set $[0,1]^{\Gamma }$.
We denote by $e_{\gamma }$ the restriction of $\pi_{\gamma }$
projection of $[0,1]^{\Gamma }$ on the space $\cals E$,
and we consider the space $X$ generated by the vectors
$\{ e_{\gamma }\}_{\gamma \in \Gamma }$ in the space $C(\cals E)$.
Then $X$ is a subspace of a W.C.G. space but it is not itself
a W.C.G. space.}

\sprth
Before we prove the theorem we give a combinatorial lemma
which is necessary in proving the result.

\thm{Lemma}
{\em Let $\Gamma $ be a set and $\{ f_{\delta }\}_{\delta
\in \Delta }$ a family
of real-valued functions defined on $\Gamma $ and
satisfying the following conditions:
\begin{list}{(\roman{lista})  }{\usecounter{lista}}
\item For every $\delta $ in $\Delta $, the function
 $f_{\delta }$ has countable support.
\item For every $\gamma \in \Gamma $ the set
$\{\delta :f _{\delta }(\gamma )\neq 0\}$ is non empty and countable.
\end{list}

Then:\\
There exist partitions $(\Gamma _{d})_{d\in D}, (\Delta _{d})_{d\in D}$
of $\Gamma $ and $\Delta $, respectively, \st $\Gamma _{d}, \Delta _{d}$
are countable and if $d_{1}\neq d_{2}, \Gamma _{d_{1}}\cap \Gamma _{d_{2}}
=\emptyset$, $\Delta _{d_{1}}\cap \Delta _{d_{2}}=\emptyset$ and
$f_{\delta }(\gamma )=0$ for $\gamma \in \Gamma _{d_{i}},
\delta \in \Delta _{d_{j}}$,
where $i=1,j=2$ or vice-versa.}

\proof
We shall show that \tx $\Gamma _{1},\Delta _{1}$ countable subsets of
$\Gamma $ and $\Delta $, respectively, \st $\Gamma _{1}=
\cup\{ {\rm supp}f_{\delta }:\delta \in \Delta _{1}\}$ and
$\Delta _{1}=\{ \delta \in \Delta :{\rm there\; exists}\; \gamma
\in \Gamma _{1}, \; f_{\delta }(\gamma )\neq 0\}$. If this has been done,
then inductively we produce the desired families
$(\Gamma _{d})_{d\in D}, (\Delta _{d})_{d\in D}$. To find
the sets $\Gamma _{1}, \Delta _{1}$ we use the standard saturation
argument. Indeed, we start with any $\delta _{0}\in \Delta $.
Set $\Delta ^{0}=\{ \delta _{0}\}$ and $\Gamma ^{0}=\{ \gamma \in \Gamma :
f_{\delta _{0}}(\gamma )\neq 0\}$. Then $\Gamma ^{0}$ is countable
and the set $\Delta ^{1}=\{ \delta \in \Delta :{\rm there\; exists}\;
\gamma \in \Gamma ^{0}\;{\rm with}\; f_{\delta }(\gamma )\neq 0\}$ is
also countable. Next we define
$\Gamma ^{1}=\{ \gamma \in \Gamma :{\rm there\; exists}\;
\delta \in \Delta ^{1}\;{\rm with}\; f_{\delta }(\gamma )\neq 0\}$.
We thus inductively produce
$\Delta ^{0}\sbs \Delta ^{1}\sbs\cdots\sbs \Delta ^{k}\sbs\cdots$,
$\Gamma ^{0}\sbs \Gamma ^{1}\sbs\cdots\sbs \Gamma ^{k}\sbs\cdots$ \st
if $\delta \in \Delta ^{k}$, the set $\{\gamma \in
\Gamma :f_{\delta }(\gamma )
\neq 0\}\sbs \Gamma ^{k}$ and if
$\gamma \in \Gamma ^{k}$, the set $\{\delta \in
 \Gamma :f_{\delta }(\gamma )
\neq 0\}\sbs \Delta ^{k+1}$. It is clear that the sets
$\Gamma _{1}=\cup_{k=1}^{\infty}\Gamma ^{k}$ and $\Delta _{1}=
\cup_{k=1}^{\infty}
\Delta ^{k}$ are the desired.\hfill Q.E.D.

\sprth {\bf Proof of the Theorem.}
Since $\cals E$ is an Eberlein compact set, $C(\cals E)$
is a W.C.G. space [A-L] and hence $X$ is, indeed, a
subspace of a W.C.G. space.

To see that $X$ is not a W.C.G. space, notice first
that for $s$ in \cals , $s\in\cals _{n}$ for some $n\in {\bf N}$
hence the functional $\frac{1}{n}s^{*}$ defined by
\[ \frac{1}{n}s^{*}\left(\sumin \alpha _{i}e_{\gamma _{i}}\right)
=\frac{1}{n}\sum\{ \alpha _{i}:\gamma _{i}\in s\}\]
is continuous and has $\ \frac{1}{n}s^{*}\\leq 1$.
Therefore $\ s^{*}\\leq n$. In particular,
since $s_{\gamma }=\{ \gamma \}\in \cals _{1}$, the functional
$s_{\gamma }^{*}$ has $\ s^{*}_{\gamma }\\leq 1$
for all $\gamma \in \Gamma $.
It is easy to see that the family $\{ s^{*}_{\gamma }\}
_{\gamma \in \Gamma }$ is $w^{*}$-total, $w^{*}$-discrete and
$\{ s^{*}_{\gamma }\}_{\gamma \in \Gamma }\cup\{ 0\}$ is $w^{*}$-compact.

Assume now that $X$ is a W.C.G. space. Then \txs a total
 $\{ y_{\delta }\}_{\delta \in \Delta }$ subset of $X$ which is weakly
discrete, and $\{ y_{\delta }\}_{\delta \in \Delta }\cup\{ 0\}$ is
weakly compact. Hence every countable subset of it defines
a weakly null \sq in the space $X$. We apply Lemma 3.7 to
the families $\{ y_{\delta }\}_{\delta \in \Delta }$ and
$\{ s^{*}_{\gamma }\}_{\gamma \in \Gamma }$ to find partitions
$\{ \Delta _{d}\}_{d\in D},
\{ \Gamma _{d}\}_{d\in D}$ of $\Delta $ and
$\Gamma $, respectively, satisfying the conclusions of the lemma.

Each $\Gamma _{d}$ is countable hence we enumerate it as
$\Gamma _{d}=\{ \gamma ^{n}_{d}\} _{n\in{\bf N}}$. We define a partition of
$\Gamma $ into a countable family $\{ \Gamma _{(n,k)} \}_{
(n,k)\in {\bf N}\times {\bf N}}$ so that $\gamma \in \Gamma _{(n,k)}$
 if and only if \txs $d\in D$ with $\gamma ^{n}_{d}=\gamma $ and \txs
$\delta\in \Delta _{d}$ with $ e^{*}_{\gamma }(y_{\delta })\geq
\frac{1}{k}$.

Since \cals\ is a quasi-Eberlein set, \txs $s\in\cals$ and
 $(n_{0},k_{0})\in {\bf N}\times {\bf N}$
\st $s\cap \Gamma _{(n_{0},k_{0})}$ is
an infinite set and $\# s\cap \Gamma _{d}\leq 1$ for all
$d\in D$ (Definition 3.1 (iii)). For every $\gamma \in
s\cap G_{(n_{0},k_{0})}$ \txs unique $d\in D$ \st
$\gamma =\gamma ^{n_{0}}_{d}$ and \txs at least one $\delta _{\gamma }
\in \Delta _{d}$ \st $ e^{*}_{\gamma }(y_{\delta _{\gamma }})\geq
\frac{1}{k}$. Consider the countable set
$\{ y_{\delta _{\gamma }}:\gamma \in s\cap \Gamma _{(n_{0},k_{0})}\}$.
As we have noticed, this set defines a weakly null
sequence. On the other hand, $s^{*}(y_{\delta _{\gamma }})
=s^{*}_{\gamma }(y_{\delta _{\gamma }})\geq\frac{1}{k}$
(by condition (i) of Definition 3.1)
for all $\gamma $ in $s\cap \Gamma _{(n_{0},k_{0})}$. This is a
contradiction and the proof is complete.\\

\vspace{-0.05in}
\hfill Q.E.D.

\thm{Remarks} (a) It is easy to see that if \cals\ is
an adequate family, i.e. if $t$ is a subset of $s$
and $s$ is in \cals\  then $t$ is also in \cals ,
then the family $\{ e_{\gamma }\}_{\gamma \in \Gamma }$ that
generates the space $X$ above is an unconditional basis.
Hence in the case of Talagrand's space, the space $X$ has
an unconditional basis.

(b) It is well known that every subspace of a W.C.G. Banach
space is a weakly \calk -analytic space [T]. Rosenthal
has proved that under Martin's Axiom  every subspace $X$ of
$\rml ^{1}(\mu )$, where $\mu$ is a probability measure,
with ${\rm dim}X< 2^{\omega }$ is a W.C.G. space. Recently S.
Merkourakis proved that this is a general phenomenon in
the class of weakly \calk -analytic Banach spaces. In
particular, he showed that under ${\rm MA}$
every weakly \calk -analytic Banach space $X$ with
${\rm dim}X< 2^{\omega }$ is a W.C.G. Banach space.

\sprth
In the sequel we denote by $RE$ the Eberleinization of
Recni\v{c}enko's space as it is defined
in 3.5. We define the James-$RE$
norm.

\thm{Definition}
Let $\phi $ be a finitely supported real valued function on
$\Gamma $. For $g$ in $RE$ we denote by $<\phi ,g>$ the real
number $\sum_{a\in g}\phi (a)g(a)$. The $J-RE$ norm of the
function $\phi $ is defined as
\[ \ \phi \ =\sup\left(\sumin <\phi ,g_{i}>^{2}\right)^{1/2},\]
where $g_{1},\ldots ,g_{n}$ are in $RE$ with \pair supports.

The space $J-RE$ is the completion in $J-RE$ norm of the
linear space $\Phi$ of all finitely supported functions
defined on $\Gamma $.

\thm{Theorem}
{\em The space $J-RE$ is a subspace of a W.C.G. space, it does not contain
isomorphically $\ell ^{1}$ and it is not a W.C.G. space.}

\proof
The proof of the theorem is analogous to the proof of certain
previous results. We shall therefore only indicate how we
get the result.

{\bf (a) The space $J-RE$ does not contain $\ell ^{1}$.}

To prove this we repeat step by step the proof of Theorem 2.2.
The only change we have to make is to replace the functional
$s^{*}$ by the functional $g^{*}$, where $g\in \cals E$.

{\bf (b) The space $J-\cals E$ is a subspace of a W.C.G. space.}

As in the proof of Theorem 2.2, the space $J-\cals E$ is a
subspace of the space $C(\overline{D}^{w^{*}}),$ where in our
case the set $D$ is defined as
\beq D&=&\left\{\sumin \lambda _{i}g_{i}^{*}:
\lambda _{i}\in\reals , g_{i}\in\cals E,
\{ g_{i}\}_{i=1}^{n}\;{\rm have\; pairwise\; disjoint\;
supports}\right.\\
& &\;\;\;\;\;\;\;\;\;\;\;\;\;\;\;\;\;
{\rm and}\;\left. \sumin \lambda ^{2}_{i}\leq 1\right\} .\eeq
Then, as in Lemma 2.4, $(\overline{D}^{w^{*}},w^{*})$ is the
continuous image of $V\times E$, where $E$ is a closed
subset of $\cals E^{{\bf N}}$ and $V$ is a closed subset of
$(B_{\ell ^{2}},w)$. Hence $V\times E$ is an Eberlein compact
set and the same holds for the set $(\overline{D}^{w^{*}},
w^{*})$ [B-R-W].
Hence $C(\overline{D}^{w^{*}})$ is W.C.G. [A-L] and the proof
is complete.

{\bf (c) The space $J-\cals E$ is not a W.C.G. space.}

For $s$ in $R$ \txs $n\in {\bf N}$ \st $s\in R_{n}$,
hence $\ g^{*}\ = \frac{1}{n}s^{*} \leq 1$.
Therefore $s^{*}\in J-RE^{*}$ and $\ s^{*}\\leq n$.
Also, for $\gamma $ in $\Gamma $, $\gamma $ is in $R$, hence if
$s^{*}_{\gamma }=\{ \gamma \} ^{*}$, we have that
$\ s^{*}_{\gamma }\\leq 1$. It is easy to see that
$\{ s^{*}_{\gamma }\}_{\gamma \in \Gamma }$ is $w^{*}$-total in $J-RE^{*}$,
$w^{*}$-discrete and $\{ s^{*}_{\gamma }\}_{\gamma \in \Gamma }\cup\{ 0\}$
is $w^{*}$-compact. Assuming now that $J-RE$ is a W.C.G.
space, \txs a subset $\{ y_{\delta }\}_{\delta \in \Delta }$
of it weakly total, weakly discrete and
$\{ y_{\delta }\}_{\delta \in \Delta }\cup\{ 0\}$ weakly compact.

What remains is to repeat exactly in the same manner the
part of Theorem 3.6 showing that
$\{ y_{\delta }\}_{\delta \in \Delta }$ has a countable
subset which does not define a weakly null sequence, thus
reaching a contradiction and completing the proof.\hfill Q.E.D.

\[{\bf REFERENCES }\]
\small
\begin{description}
\parsep =0pt
\itemsep =0pt
\parskip=0pt
\bb  {[A-L]}]
  {\sc D.Amir} and {\sc J. Lindenstauss,}
  {\em The structure of weakly compact sets in Banach spaces,}
   Ann. of Math. 88(1968) 35-46.
\bb {[A-M]}]
  {\sc S. Argyros} and {\sc S. Merkourakis,}  {\em On weakly Lindel\"{o}f
Banach spaces} Rocky Mountain J. of Math. 23 (1993) 395-446.
\bb {[A-M-N]}] {\sc S. Argyros, S. Merkourakis} and {\sc S. Negrepontis}
{\em Functional-analytic properties of Corson-compact spaces,}
Studia Math. 89(1988) 197-229.
\bb {[A-N]}] {\sc S. Argyros} and {\sc S. Negrepontis,}
{\em On weakly \calk -countably determined spaces of continuous
functions,} Proc. A.M.S. 87 (1983) 731-736.
\bb {[B-R-W]}] {\sc Y. Benyamini, M.E. Rudin} and {\sc M. Wage,}
{\em Continuous Images of weakly compact sets in Banach spaces,}
Pacific J. of Math. 70 (1977) 309-324.
\bb {[D-G-Z]}] {\sc R. Deville, G. Godefroy} and {\sc V. Zizlev,}
{\em Smoothness and renormings in Banach spaces}
Longman (1992).
\bb {[F]}] {\sc M. Fabian,} {\em Each weakly countably determined
Asplund space admits a Frechet differentiable norm,}
Bull. Austral. Math. Soc. 36  (1987) 367-374.
\bb {[Fa]}] {\sc V. Farmaki,} {\em The structure of Eberlein
uniformly Eberlein and Talagrand compact spaces in
$\sum (\reals^{\Gamma} )$.}
Fundamenta Math. 128 (1987) 15-28.
\bb {[J]}] {\sc R.C. James,} {\em A seperable somewhat reflexive
Banach space with non-seperable dual,} Bull. AMS 80 (1974)
738-743.
\bb {[M-N]}] {\sc S. Merkourakis} and {\sc S. Negrepontis,}
{\em Banach spaces and Topology II,} In M. Hu\v{s}ek and
J. van Mill (1992) 495-536.
\bb {[Ny]}] {\sc P. Nyikos,} {\em Classes of compact sequential
spaces,} Lecture Notes in Math. 1401 Springer Verlag (1988)
135-159.
\bb {[O-S-V]}] {\sc J. Orihuela, W. Schachermayer} and {\sc M. Valdivia,}
{\em Every Radon-Nikodym Corson compact space is Eberlein compact,}
Studia Math. 98 (1991) 157-174.
\bb {[Rai]}] {\sc J. Rainwater,} {\em Weak convergence of bounded sequences,}
 Proc. A.M.S. 14 (1963) 999.
\bb {[R]}] {\sc H.P. Rosenthal,} {\em The heredity problem for
weakly compactly generated Banach spaces,}
Compos. Math. 28 (1974) 83-111.
\bb {[So]}] {\sc G.A. Sokolov,} {\em On some classes of compact
spaces lying in $\sum$-products,}
Comm. Math. Univ. Carolinae 25 (1984) 219-231.
\bb {[S]}] {\sc C. Stegall,} {\em More facts about conjugate Banach
spaces with the Radon-Nikodym property,} Acta Univ. Carol. Math.
Phys., 31 (1990) 107-117.
\bb {[T]}] {\sc M. Talagrand,} {\em Espaces de Banach faiblement
K-analytiques,} Ann. of Math. 110 (1979) 407-438.
\bb {[To]}] {\sc S. Todorcevi\v{c},} {\em Trees and linearly
ordered sets,} In Kunen and Vaughan (1984) 235-293.
\bb {[V]}] {\sc M. Valdivia,} {\em Resolutions of the identity in certain
Banach spaces,} Collect. Math. 39 (1988) 127-140.

\end{description}
\hfill  DEPARTMENT OF MATHEMATICS

\vspace{-0.05in}
\hfill  UNIVERSITY OF ATHENS\

\vspace{-0.05in}
\hfill
 157 84 Athens GREECE

\end{document}